\documentclass{sig-alternate}

\newtheorem{lemma}{Lemma}
\newtheorem{theorem}[lemma]{Theorem}
\newtheorem{corollary}[lemma]{Corollary}
\newtheorem{proposition}[lemma]{Proposition}
\newtheorem{definition}{Definition}

\def\ba{\begin{array}}
\def\ea{\end{array}}
\def\bi{\begin{itemize}}
\def\ei{\end{itemize}}
\def\bd{\begin{description}}
\def\ed{\end{description}}
\def\bu{\begin{enumerate}}
\def\eu{\end{enumerate}}
\def\be{\begin{equation}}
\def\ee{\end{equation}}
\def\bn{\begin{eqnarray}}
\def\en{\end{eqnarray}}

\def\bproof{\noindent{\bfseries\itshape Proof.}\hskip1pc\ignorespaces}
\def\endproof{{\hfill \Large $\bf\bigstrut \Box$}\vskip .15cm}

\def\bt{\begin{theorem}\rm}
\def\et{\end{theorem}}
\def\bc{\begin{corollary}\rm}
\def\ec{\end{corollary}}
\def\bl{\begin{lemma}\rm}
\def\el{\end{lemma}}
\def\bp{\begin{proposition}\rm}
\def\ep{\end{proposition}}
\def\bdf{\begin{definition}\rm}
\def\edf{\end{definition}}
\def\bx{\begin{example}}
\def\ex{\end{example}}

\def\qed{\hbox{\hfill []}}
\def\cal{\mathcal}

\def\P{{\bf P}}
\def\Q{{\bf Q}}

\def\S{{\bf S}}

\def\V{{\bf V}}

\def\e{{\bf e}}

\def\i{{\bf i}}
\def\j{{\bf j}}
\def\k{{\bf k}}

\def\p{{\bf p}}
\def\q{{\bf q}}

\def\s{{\bf s}}

\def\v{{\bf v}}

\def\x{{\bf x}}

\def\0{{\bf 0}}
\def\1{{\bf 1}}
\def\2{{\bf 2}}
\def\3{{\bf 3}}
\def\4{{\bf 4}}
\def\5{{\bf 5}}
\def\6{{\bf 6}}
\def\7{{\bf 7}}
\def\8{{\bf 8}}
\def\9{{\bf 9}}

\def\bq{\bar{\bf q}}
\def\barp{\bar{\bf p}}
\def\bQ{\bar{\bf Q}}

\def\rG{\hbox{\rm G}}

\newbox\bigstrutbox
\setbox\bigstrutbox=\hbox{\vrule height12.5pt depth5pt width0pt}
\def\bigstrut{\relax\ifmmode\copy\bigstrutbox\else\unhcopy\bigstrutbox\fi}

\newbox\Bigstrutbox
\setbox\Bigstrutbox=\hbox{\vrule height17.5pt depth5pt width0pt}
\def\Bigstrut{\relax\ifmmode\copy\Bigstrutbox\else\unhcopy\Bigstrutbox\fi}

\setpagenumber{1}

\begin{document}

\title{Normalization of Quaternionic Polynomials}

\numberofauthors{1}

\author{\alignauthor Hongbo Li,\ \ \ \ Lei Huang,\ \ \ \ Yue Liu\\ 
\vskip .1cm
       \affaddr{KLMM, AMSS, Chinese Academy of Sciences, Beijing 100190, China}\\
\vskip .1cm
       \email{\small
       hli@mmrc.iss.ac.cn, lhuang@mmrc.iss.ac.cn, liuyue10@mails.gucas.ac.cn}   
}

\maketitle

\begin{abstract} \vskip .1cm
Quaternionic polynomials are generated by quaternionic variables and the quaternionic product.
This paper proposes the generating ideal of quaternionic polynomials in tensor algebra, 
finds the Gr\"obner base of the ideal in the case of pure imaginary quaternionic 
variables, and describes the normal forms of such quaternionic polynomials explicitly.  
\end{abstract}

\category{I.1.1}{Computing Methodologies}{SYMBOLIC AND ALGEBRAIC MANIPULATION}

\terms{Theory}

\keywords{Quaternionic polynomial ring; syzygy ideal; non-commutative Gr\"obner base; normal form; 
Clifford algebra.}

\section{Introduction}
\vskip .3cm

Quaternions are one of the simplest examples of a non-commutative associative algebra. When pure imaginary qua{-}ternions 
are used to represent 3D vectors, the adjoint action of the group of unit quaternions upon them provides a 
representation of all 3D rotations.  
Because of this, quaternions have important applications in 3D Euclidean computing
\cite{altmann}. 

In symbolic geometric computing, when geometric entities and their transformations are represented by 
coordinate-free symbols, algebraic manipulations among the symbols often lead
to rather simple and geometrically interpretable results \cite{li}. For example, using $\v$ instead of 
$x\i+y\j+z\k$ as a vector variable, provides a typical coordinate-free
representation of a 3D direction (or point) by a pure imaginary quaternion. Under the 
quaternionic product, these variables generate a ring of non-commutative polynomials. 
It is natural to ask for symbolic algebraic algorithms to manipulate such {\it quaternionic polynomials}. 

When the coefficients of the quaternionic polynomials are taken in the field  of rational numbers $\mathbb Q$, 
the quaternionic product, being associative and multilinear, must be a quotient of the $\mathbb Q$-tensor product
modulo an ideal $\cal I$ whose generators are some tensors. Finding such generators and proving that they
are sufficient to define the quaternionic product is the first task to finish.
The first theorem in this paper provides very simple generators of $\cal I$, called the {\it syzygy ideal}
of quaternionic polynomials. 

The next task is naturally on computing the Gr\"{o}bner bases of $\cal I$ and then the
normal forms of quaternionic polynomials.
This can always be done using the classical theory and algorithm on non-commutative Gr\"{o}bner bases \cite{mora}.
For a fixed number $n$ of quaternionic variables, one can compute a Gr\"{o}bner base of the syzygy ideal
${\cal I}$ in the $n$ variables, and for every degree $d$, the number of degree-$d$ polynomials in the
computed Gr\"{o}bner base is finite. After investigating the elements in the Gr\"{o}bner bases for small values
of $n$ and $d$, we come up with some observations on the general form of the Gr\"{o}bner base for arbitrary $n$
and $d$, and of the normal forms of quaternionic polynomials with respect to the Gr\"{o}bner base.

The main part of this paper is the establishment of two theorems that provide explicitly
(1) a reduced Gr\"{o}bner base of the syzygy ideal ${\cal I}$ of quaternionic polynomials
in $n$ pure imaginary quaternionic variables, (2) the structure of the normal forms 
relative to the Gr\"{o}bner base. The results are extremely simple. For example, each polynomial in the
Gr\"{o}bner base has either 2 or 4 terms; 
the sequence of variables in a monomial of normal form always forms two
intertwining non-descending subsequences in terms of the ordering of the variables (Figure \ref{normal} and
Figure \ref{normal:2} of this paper).

The proofs of the theorems are very long and are divided into many cases.
We provide detailed proofs for several cases for illustration, but skip the rest by limitation of
the size of this paper. The results in this paper are not only
foundations for symbolic manipulations of quaternionic polynomials, but also helpful in
exploring general results on the Clifford product among vectors of arbitrary dimension \cite{hestenes}, 
now that the dimension-three case is solved here as a by-product.

\section{Quaternionic polynomial ring \\
generated by vectors}
\setcounter{equation}{0}

Let $\q_1, \q_2, \ldots, \q_n$ be a set of symbols representing quaternionic variables, and let 
$\bq_1, \bq_2, \ldots, \bq_n$ be the corresponding quaternionic conjugates.
A {\it vector} $\v$ refers to a pure imaginary quaternion, {\it i.e.}, $\bar{\v}=-\v$,
and a {\it scalar} $\s$ refers to a real quaternion, {\it i.e.}, $\bar{\s}=\s$.

A {\it monic quaternionic monomial} in variables
$\q_1, \q_2, \ldots, \q_n$ refers to the quaternionic product of a permutation of the elements in
a multiset of the $2n$ symbols $\{\q_i, \bq_i\,|\,i=1,\ldots,n\}$. A   
{\it quaternionic polynomial} with rational coefficients is a linear combination of 
monic quaternionic monomials with coefficients in the rational numbers field $\mathbb Q$.
We always use juxtaposition of elements to represent their quaternionic product.
For a monomial $\Q=\q_{i_1}\q_{i_2}\cdots \q_{i_k}$, its {\it reversion} is
\be
\Q^\dagger:=\q_{i_k}\cdots \q_{i_2}\q_{i_1}.
\ee




The quaternionic product is the quotient of the $\mathbb Q$-tensor product modulo some {\it syzygy ideal}. 
Since in symbolic manipulation, a scalar
is characterized by its commutativity with everything, we first consider
characterizing vectors and their quaternionic products by the syzygy ideal, 
then extend the results to general quaternionic variables.

Let $\v_1, \v_2, \ldots, \v_n$ be vector variables. They belong to a 3D real inner-product space
with metric diag$(-1,-1,-1)$. As the standard basis $\{\i,\j,\k\}$ cannot be constructed from the given
vector variables without using square-root and division operations, they are not allowed in manipulating the
quaternionic polynomials of the $\v_i$. 

Let $\v_1\prec \v_2\ldots \prec \v_n$. For a monic quaternionic monomial 
$\v_{i_1}\v_{i_2}\cdots \v_{i_k}$, the {\it leading variable} refers to $\v_{i_1}$, and the {\it trailing variable}
refers to $\v_{i_k}$. The monomial is said to be {\it non-descending} if $i_1\leq i_2\leq \ldots\leq i_k$,
and is said to be
{\it ascending} if $i_1< i_2< \ldots< i_k$. The {\it degree}, or {\it length}, of the monomial is $k$.
The {\it degree} of a quaternionic polynomial is the biggest degree of its monomials disregarding the 
$\mathbb Q$-coefficients. The {\it lexicographic ordering} among monomials is always assumed.
The {\it leading term} of a quaternionic polynomial is the term of highest order.

The {\it inner product} of two vectors (not necessarily vector variables) $\v_i, \v_j$
is defined by
\be
[\v_i\v_j]:=(\v_i\v_j+\v_j\v_i)/2. \label{inner}
\ee
The result is a scalar, so it commutes with a third vector $\v_k$:
$[\v_i\v_j]\v_k=\v_k[\v_i\v_j]$.

The {\it vector product} of two vectors $\v_i, \v_j$
is defined by
\be
\v_i\times \v_j:=(\v_i\v_j-\v_j\v_i)/2.
\ee
The result is a vector, so its inner product with a third vector $\v_k$ is a scalar. Define
\be
[\v_i\v_j\v_k]:=(\v_i\v_j\v_k-\v_k\v_j\v_i)/2.
\ee
Then $[\v_i\v_j\v_k]=[(\v_i\times \v_j)\v_k]$ is a scalar, so for a fourth vector $\v_l$,
$[\v_i\v_j\v_k]\v_l=\v_l[\v_i\v_j\v_k]$.

Define
\be
[\v_i\v_j\v_k\v_l]:=(\v_i\v_j\v_k\v_l+\v_l\v_k\v_j\v_i)/2.
\ee

\bt \label{thm:1}
Let $\v_1, \v_2, \ldots, \v_n$ be $n$ different symbols, where $n\geq 2$. Define
the product among them, denoted by juxtaposition of elements, as the 
$\mathbb Q$-tensor product modulo the 
two-sided ideal generated by the following tensors:
\be\hskip -.15cm
\ba{ll}
{\rm V2}:   &\hskip .1cm \v_i\otimes \v_i\otimes \v_j- \v_j\otimes \v_i\otimes \v_i; \\
{\rm V3}:   &\hskip .1cm (\v_i\otimes \v_j+\v_j\otimes \v_i)\otimes \v_k-\v_k\otimes (\v_i\otimes \v_j+\v_j\otimes \v_i); \\
{\rm V4}:   &\hskip .1cm (\v_i\otimes \v_j\otimes \v_k-\v_k\otimes \v_j\otimes \v_i)\otimes \v_l\\
& \hskip .8cm
-\v_l\otimes (\v_i\otimes \v_j\otimes \v_k-\v_k\otimes \v_j\otimes \v_i),
\ea
\label{eq0}
\ee
for any $i\neq j\neq k\neq l$ in $1,2,\ldots, n$. Denote by $\cal Q$ the $\mathbb Q$-algebra defined by the above product 
and generated by the $\v_i$. Denote 
\be
\overrightarrow{\mathbb Q}:={\mathbb Q}(\{\v_i\otimes \v_l+\v_l\otimes \v_i, 
\v_i\otimes \v_j\otimes \v_k-\v_k\otimes \v_j\otimes \v_i\,|\,i\neq j\neq k\}).
\label{q:field}
\ee
Then 
(1) the $\v_i$ and $\v_i\times \v_j$ span a 3D $\overrightarrow{\mathbb Q}$-vector space ${\cal V}^3$; \\
(2) if all the elements in (\ref{q:field}) extending $\mathbb Q$
are real numbers, and the inner product of real vector space
${\cal V}^3$ induced from (\ref{inner}) is definite, then
the product among the $\v_i$ is the quaternionic product, {\it i.e.},
$\cal Q$ is the quaternionic polynomial ring generated by the $\v_i$ as vector variables.
\et

\bproof The proof is composed of the following steps:

{\it Step 1}. For all $1\leq i,j,k,l\leq n$, 
\be\ba{rll}
\v_i\v_i\v_j &=& \v_j\v_i\v_i; \\
(\v_i\v_j+\v_j\v_i)\v_k &=& \v_k(\v_i\v_j+\v_j\v_i); \\
(\v_i\v_j\v_k-\v_k\v_j\v_i)\v_l &=& \v_l(\v_i\v_j\v_k-\v_k\v_j\v_i).
\ea
\label{eq1}
\ee
The proof is easily established by verifying that V3 for $i=k$, and V4 for $l=i$ or $j$ are in the ideal generated by 
elements in (\ref{eq0}). So the center of the ring $\cal Q$ is
the subring generated by the $[\v_i\v_i]$, $[\v_i\v_j]$, $[\v_i\v_j\v_k]$ for all $1\leq i<j<k\leq n$.

Obviously, the latter two equalities in (\ref{eq1}) can be written as the following {\it shift invariance} of
$[\v_i\v_j\v_k]$ and $[\v_i\v_j\v_k\v_l]$:
\be
[\v_i\v_j\v_k] = [\v_k\v_i\v_j], \hskip .5cm
[\v_i\v_j\v_k\v_l] = [\v_l\v_i\v_j\v_k].
\label{eq2}
\ee

{\it Step 2}. Since 
$\v_i\otimes \v_i$, $\v_i\otimes \v_j+\v_j\otimes \v_i$ and 
$\v_i\otimes \v_j\otimes \v_k-\v_k\otimes \v_j\otimes \v_i$ for all $1\leq i<j<k\leq n$ are not in the ideal
generated by elements V2, V3, V4,
$[\v_i\v_i]$, $[\v_i\v_j]$, $[\v_i\v_j\v_k]$ for all $1\leq i<j<k\leq n$ are nonzero.

Consider the following {\it Cramer's rule}: for all $i,j,k,l$,
\be
[\v_i\v_j\v_k]\v_l-[\v_i\v_j\v_l]\v_k+[\v_i\v_k\v_l]\v_j-\v_i[\v_j\v_k\v_l]=0. \label{cramer:1}
\ee
By the first equality in (\ref{eq2}), the left side of (\ref{cramer:1})
can be simplified to $[\v_i\v_j\v_k]\v_l-\v_l[\v_i\v_j\v_k]$,
which equals zero by the last equality in (\ref{eq1}).
So if $n>2$, the $\v_i$ span a 3D $\overrightarrow{\mathbb Q}$-vector space. 

{\it Step 3}. The following relations can be verified using (\ref{eq1}):
\be\ba{rcl}
(\v_i\times \v_j)\times \v_k &=& [\v_j\v_k]\v_i-[\v_i\v_k]\v_j;\\
{[}(\v_i\times \v_j)(\v_k\times \v_l)] &=& [\v_i\v_l][\v_j\v_k]-[\v_i\v_k][\v_j\v_l];\\
(\v_i\times \v_j)\times (\v_k\times \v_l) &=& [\v_j\v_k\v_l]\v_i-[\v_i\v_k\v_l]\v_j.
\ea
\label{eq3}
\ee
Notice that in the corresponding formulas for the vector algebra over ${\mathbb R}^3$, the signs on the right side
of (\ref{eq3}) are opposite. In fact, the expressions in the three lines of (\ref{eq3}) are respectively equal to
\[
\frac{\v_i\v_j\v_k-\v_k\v_i\v_j}{2}; \hskip .2cm 
[\v_i\v_j(\v_k\times \v_l)]; \hskip .2cm
\frac{\v_i\v_j\v_k\v_l-\v_k\v_l\v_i\v_j}{2}.
\]

If $n=2$, it is easy to see that $\v_1\otimes \v_2-\v_2\otimes \v_1$, and 
$\v_1\otimes \v_2\otimes \v_1\otimes \v_2
+\v_2\otimes \v_1\otimes \v_2\otimes \v_1-2 \v_1\otimes \v_1\otimes \v_2\otimes \v_2$
are not in the defining ideal of $\cal Q$, so $\v_1\times \v_2\neq 0$, and
\be\hskip -.1cm
\ba{lll}
{[}(\v_1\times \v_2)\v_1\v_2] &\hskip -.2cm=&\hskip -.2cm (\v_1\times \v_2)(\v_1\times \v_2)\\
&\hskip -.2cm
=&\hskip -.2cm
 (\v_1\v_2\v_1\v_2+\v_2\v_1\v_2\v_1-2\v_1\v_1\v_2\v_2)/4\\
&\hskip -.2cm\neq &\hskip -.2cm 0.
\ea
\label{n2}
\ee 
Then $\v_1, \v_2, \v_1\times \v_2$ are linearly independent over $\overrightarrow{\mathbb Q}$.

By (\ref{eq1}) and (\ref{eq3}), the following {\it Cramer's rule} can be verified: 
for all $1\leq i,j,k,l,m\leq n$,
\be\ba{lll}
\hskip -.2cm
[(\v_i\times \v_j)\v_k\v_l]\v_m-[(\v_i\times \v_j)\v_k\v_m]\v_l\\
\hskip .1cm
+[(\v_i\times \v_j)\v_l\v_m]\v_k
-[\v_k\v_l\v_m](\v_i\times \v_j) &\hskip -.2cm=&\hskip -.2cm 0.
\ea
\label{cramer:2}
\ee
Hence the $\v_i$ and $\v_j\times \v_k$ span a 3D $\overrightarrow{\mathbb Q}$-vector space.

{\it Step 4}. The inner product defined by (\ref{inner})
can be extended bilinearly to all elements in ${\cal V}^3$. Then $\v_1\times \v_2$ is orthogonal to both
$\v_1$ and $\v_2$. Furthermore, the inner products of the three elements with themselves separately are all nonzero.
When the $[\v_i\v_i]$, $[\v_i\v_j]$, $[\v_i\v_j\v_k]$ are all real numbers, let 
\[
\e_1=\v_1/\sqrt{|[\v_1\v_1]|}, \hskip .24cm
\e_2=\v_1\times \v_2/\sqrt{|[(\v_1\times \v_2)(\v_1\times \v_2)]|},
\]
and let $\e_3=\e_1\times \e_2$. Then $\{\e_1, \e_2, \e_3\}$ form an orthogonal basis of ${\cal V}^3$, and
$[\e_i\e_i]=\pm 1$.
Furthermore, $\e_3^2=[\e_3\e_3]=-\e_1^2\e_2^2$. So the metric of ${\cal V}^3$ is either 
diag$(-1,-1,-1)$ or diag$(-1,1,1)$.

{\it Step 5}. It is easy to prove $\e_2\e_3=\e_2\times \e_3=-\e_2^2\e_1$ and $\e_3\e_1=\e_3\times \e_1=-\e_1^2\e_2$. 
Moreover, $\e_1\e_2=\e_1\times \e_2=\e_3$. So
if $\e_i^2=-1$ for $i=1,2,3$, then the product in ${\cal V}^3$ 
is the quaternionic product.
\endproof

As a consequence, all polynomial identities among the $\v_i$ in the ring $\cal Q$ are generated by 
relations of the form $[\v_i\v_j]\v_k=\v_k[\v_i\v_j]$ and $[\v_i\v_j\v_k]\v_l=\v_l[\v_i\v_j\v_k]$, or equivalently,
by relations of the form $[\v_i\v_j\v_k]=[\v_k\v_i\v_j]$ and $[\v_i\v_j\v_k\v_l]=[\v_l\v_i\v_j\v_k]$.
Neither the assumption that the $[\v_i\v_i]$, $[\v_i\v_j]$ and $[\v_i\v_j\v_k]$ are real numbers, nor the assumption
that the metric of ${\cal V}^3$ is definite, can be represented by polynomial equalities of the $\v_i$. 

Henceforth we refer to the ring $\cal Q$ defined in Theorem \ref{thm:1} as the 
{\it quaternionic polynomial ring} generated by vector variables $\v_i$, and call the defining ideal of 
$\cal Q$ the {\it syzygy ideal}, and denote it by 
${\cal I}[\v_1, \v_2, \ldots, \v_n]$. 

Let $\V_k=\v_{j_1}\v_{j_2}\cdots \v_{j_k}$. Define
\be
[\V_k]:=(\V_k+(-1)^k \V_k)/2,\hskip .4cm
A(\V_k):=(\V_k-(-1)^k \V_k)/2.
\ee

\bp \label{prop:1}
$\v_i[\V_k]-[\V_k]\v_i$ is in ${\cal I}[\v_1, \v_2, \ldots, \v_n]$. 
The commutation relation $\v_i[\V_k]=[\V_k]\v_i$ can also be written as
the following shift-invariance: 
\be
[\v_i\V_k]=[\V_k\v_i]. 
\ee
The $A(\V_k)$ are in the 3D $\overrightarrow{\mathbb Q}$-vector space
spanned by the $\v_i$ and $\v_j\times \v_k$.
\ep

\bproof
The following identities can be established by induction, although not very easy: for all $l>1$,
\be\hskip -.2cm
\ba{lll}
{[}\v_{j_1}\v_{j_2}\cdots \v_{j_{2l}}] & \hskip -.2cm=& \hskip -.2cm\sum_{i=2}^{2l} 
(-1)^i [\v_{j_1}\v_{j_i}][\v_{j_2}\v_{j_3}\cdots \check{\v}_{j_i}\cdots \v_{j_{2l}}];\\

A(\V_{2l-1}) &\hskip -.2cm=& \hskip -.2cm \sum_{(2l-2,1)\vdash \V_{2l-1}} [{\V_{2l+1}}_{(1)}] {\V_{2l-1}}_{(2)}; \\

A(\V_{2l}) &\hskip -.2cm=& \hskip -.2cm \sum_{(2l-2,2)\vdash \V_{2l}} [{\V_{2l}}_{(1)}]A({\V_{2l}}_{(2)});\\

{[}\V_{2l+1}] &\hskip -.2cm=& \hskip -.2cm \sum_{(2l-2,3)\vdash \V_{2l+1}} [{\V_{2l+1}}_{(1)}][{\V_{2l+1}}_{(2)}],
\ea
\label{eq4}
\ee
where (i) $\check{\v}_{j_i}$ denotes that $\v_{j_i}$ does not occur in the sequence,
(ii) $(h,m-h)\vdash \V_{m}$ is a bipartition of the $m$ elements in the sequence $\V_m$ into two subsequences
${\V_{m}}_{(1)}$ and ${\V_{m}}_{(2)}$ of length $h$ and $m-h$ respectively, (iii) in $[{\V_{m}}_{(1)}]$, the product of
the $h$ elements in the subsequence is denoted by the same symbol ${\V_{m}}_{(1)}$, (iv) the summation
$\sum_{(h,m-h)\vdash \V_{m}}$ is over all such bipartitions of $\V_m$, and the sign of permutation of
the new sequence ${\V_{m}}_{(1)}, {\V_{m}}_{(2)}$ is assumed to be carried by the first factor $[{\V_{m}}_{(1)}]$ 
of the addend. All the conclusions of the proposition follow (\ref{eq4}) and (\ref{eq1}).
\endproof

Alternatively, we can use existing results on Clifford algebras to deduce (\ref{eq4}) directly. 
In the proof of Theorem \ref{thm:1}, we have shown that the $\v_i$ and $\v_j\times \v_k$ span
a 3D $\overrightarrow{\mathbb Q}$-vector space. To show the difference between the quaternionic product and the
Clifford product among the vectors $\v_i$ we need to reduce the field $\overrightarrow{\mathbb Q}$.

\bl \label{lemma1}
Let $\cal Q$ be the $\mathbb Q$-algebra defined by (\ref{eq0}). Then for vectors $\v_{i_1}, \ldots, \v_{i_6}$,
\be
[\v_{i_1}\v_{i_2}\v_{i_3}][\v_{i_4}\v_{i_5}\v_{i_6}]
=-\left|\ba{lll}
{[}\v_{i_1}\v_{i_4}] & [\v_{i_1}\v_{i_5}] & [\v_{i_1}\v_{i_6}]\\
{[}\v_{i_2}\v_{i_4}] & [\v_{i_2}\v_{i_5}] & [\v_{i_2}\v_{i_6}]\\
{[}\v_{i_3}\v_{i_4}] & [\v_{i_3}\v_{i_5}] & [\v_{i_3}\v_{i_6}]
\ea
\right|.
\label{lem1:eqn}
\ee
\el

\bproof
Using the second equality of (\ref{eq3}) in the reverse direction, we can write
the right side of (\ref{lem1:eqn}) after expanding the determinant as
\[\ba{ll}
& -[\v_{i_1}\v_{i_4}][(\v_{i_2}\times\v_{i_3})(\v_{i_5}\times\v_{i_6})]\\
&
+[\v_{i_2}\v_{i_4}][(\v_{i_1}\times\v_{i_3})(\v_{i_5}\times\v_{i_6})]\\
&
-[\v_{i_3}\v_{i_4}][(\v_{i_1}\times\v_{i_2})(\v_{i_5}\times\v_{i_6})]\\

=& -[\{[\v_{i_2}\v_{i_3}(\v_{i_5}\times\v_{i_6})]\v_{i_1}
-[\v_{i_1}\v_{i_3}(\v_{i_5}\times\v_{i_6})]\v_{i_2}\\
&\hfill
+[\v_{i_1}\v_{i_2}(\v_{i_5}\times\v_{i_6})]\v_{i_3}\}\v_{i_4}].
\ea
\]
By (\ref{cramer:2}), the above expression can be simplified to 
$-[\v_{i_1}\v_{i_2}$ $\v_{i_3}][(\v_{i_5}\times \v_{i_6})\v_{i_4}]$,
which equals the left side of 
(\ref{lem1:eqn}).
\endproof

Let
\be
\hat{\mathbb Q}:={\mathbb Q}(\{\v_i\otimes \v_l+\v_l\otimes \v_i, \,|\,1\leq i,l\leq n\}).
\label{q:field2}
\ee
Then when $n=2$, by (\ref{n2}), $\v_1, \v_2, \v_1\times \v_2$ span a 3D inner-product space over $\hat{\mathbb Q}$.
When $n>2$, let
\be
\iota:=[\v_1\v_2\v_3].
\ee
Then $\iota\neq 0$, and $[\v_i\v_j\v_k]=\iota\{(\iota[\v_i\v_j\v_k])/(\iota\iota)\}\in \hat{\mathbb Q}\iota$.
The center of $\cal Q$ as a $\hat{\mathbb Q}$-algebra is $\hat{\mathbb Q}+\hat{\mathbb Q}\iota$.

By multiplying (\ref{cramer:1}) with $\iota$, we get that the $\v_i$ span a 3D $\hat{\mathbb Q}$-vector space 
$\hat{\cal V}^3$. Furthermore, from (\ref{cramer:2}) we get that the $\iota \v_i\times \v_j$ are in  
$\hat{\cal V}^3$. The inner product between the $\v_i$ can be extended bilinearly to an inner product in
$\hat{\cal V}^3$. 

In the $\hat{\mathbb Q}$-Clifford algebra $Cl(\hat{\cal V}^3)$ over $\hat{\cal V}^3$,
By Theorem 5.64, Theorem 5.65 and Proposition 6.22 in \cite{li}, the two sides of the four equalities
in (\ref{eq4}) 
are different expressions of the grade-0,1,2,3 part of $\V_{2l}$, $\V_{2l-1}$, 
$\V_{2l}$, $\V_{2l+1}$ respectively. 
In $Cl(\hat{\cal V}^3)$, $\iota$ is not a scalar, but a grade-3 element. 
In the case of quaternions, the quaternionic product is a homomorphic image of the Clifford product because of
V2, V3, so (\ref{eq4}) is automatically satisfied.

The right side of the third equality in (\ref{eq4})
when multiplied by $\iota$, becomes a vector in $\hat{\cal V}^3$, while
the right side of the last equality when multiplied by $\iota$, becomes an element in $\hat{\mathbb Q}$.
These discussions lead to the following conclusion:

\bp
In the $\mathbb Q$-algebra $\cal Q$ defined by (\ref{eq0}), 
the $\v_i$ and $\iota\v_i\times \v_j$ span a 3D $\hat{\mathbb Q}$-vector space $\hat{\cal V}^3$. 
For any $k\geq 1$, $A(\v_{j_1}\v_{j_2}\cdots \v_{j_{2k+1}})$ 
and $\iota A(\v_{j_1}\v_{j_2}\cdots \v_{j_{2k}})$
are both in $\hat{\cal V}^3$.
\ep

{\it Remark.} In Clifford algebra, by Proposition 6.23 of \cite{li},
the outer product (or exterior product) of four vectors
$\v_i, \v_j, \v_k$, $\v_l$ equals one quarter of V4. So when the field is extended to
$\hat{\mathbb Q}$, V4 states that the dimension of the $\hat{\mathbb Q}$-vector space spanned by the $\v_m$
is at most three. From this aspect, V2, V3, V4 also generate the {\it syzygy ideal} of
$Cl(\hat{\cal V}^3)$, if the latter is taken as a $\mathbb Q$-algebra.

\section{Quaternionic polynomial ring}
\setcounter{equation}{0}

Let $\q_1, \ldots, \q_n$ be quaternionic variables, and let $\bq_1, \ldots, \bq_n$ be their quaternionic conjugates.
To simplify notation, we introduce $\p_i$ to represent any of $\q_i, \bq_i$, and define $\bar{\bq}_i:=\q_i$.
For a monomial $\Q=\p_{i_1}\cdots \p_{i_k}$, its {\it conjugate} is defined by
\be
\bQ:=\barp_{i_k}\cdots \barp_{i_2}\barp_{i_1}.
\ee

Define the {\it scalar part} $[\Q]$ and the {\it vector part} $A(\Q)$ of $\Q$ as following:
\be
[\Q]:=(\Q+\bQ)/2, \hskip .4cm A(\Q):=(\Q-\bQ)/2.
\ee
Then $[\q_i]$ commutes with $\p_j$, and the quaternionic product among 
vectors $\v_i:=A(\q_i)$ is characterized by (\ref{eq0}). From Theorem \ref{thm:1} and Proposition \ref{prop:1}
we can deduce the following:

\bt \label{thm:2}
Let $\q_1, \q_2, \ldots, \q_n$ and $\bq_1, \bq_2, \ldots, \bq_n$
be $2n$ different symbols, where $n\geq 2$. 
Let $\p_i$ be one of $\q_i, \bq_i$, and let $\bar{\bq}_i=\q_i$.
Define
the product among the $2n$ symbols, denoted by juxtaposition of elements, as the
$\mathbb Q$-tensor product modulo the 
two-sided ideal generated by the following tensors:
\be\hskip -.15cm
\ba{ll}
{\rm Q0}:   & \hskip -.15cm\q_i\otimes \bq_i-\bq_i\otimes \q_i; \\
{\rm Q1}:   & \hskip -.15cm(\q_i+\bq_i)\otimes \p_j - \p_j\otimes (\q_i+\bq_i); \\
{\rm Q2}:   & \hskip -.15cm
\q_i\otimes \bq_i\otimes \p_j- \p_j\otimes \q_i\otimes \bq_i; \\
{\rm Q3}:   & \hskip -.15cm(\p_i\otimes \p_j+\barp_j\otimes \barp_i)\otimes \p_k
-\p_k\otimes (\p_i\otimes \p_j+\barp_j\otimes \barp_i); \\
{\rm Q4}:   & \hskip -.15cm(\p_i\otimes \p_j\otimes \p_k+\barp_k\otimes \barp_j\otimes \barp_i)\otimes \p_l
\\
& \hskip .4cm
-\p_l\otimes (\p_i\otimes \p_j\otimes \p_k+\barp_k\otimes \barp_j\otimes \barp_i),
\ea
\label{eq5}
\ee
for any $i\neq j\neq k\neq l$ in $1,2,\ldots, n$. Denote by $\tilde{\cal Q}$ the ring defined by the above product 
and generated by the $\p_i$, and denote 
\be\hskip -.24cm
\ba{r}
\tilde{\mathbb Q}:={\mathbb Q}(\{\q_i+\bq_i, \q_i\otimes \bq_i, \p_i\otimes \p_j+\barp_j\otimes \barp_i, \hskip .8cm
\hbox{  }\\
\p_i\otimes \p_j\otimes \p_k+\barp_k\otimes \barp_j\otimes \barp_i\,|\,i\neq j\neq k\}).
\ea
\label{q:field:2}
\ee
Then 
(1) the $A(\q_i)$ and $A(\q_i)\times A(\q_j)$ span a 3D $\tilde{\mathbb Q}$-vector space $\tilde{\cal V}^3$; \\
(2) if all the elements in (\ref{q:field:2}) extending $\mathbb Q$   
are real numbers, and the inner product of real vector space
$\tilde{\cal V}^3$ induced from (\ref{inner}) for the $\v_i=A(\q_i)$ is definite, then
the product among the $\p_i$ is the quaternionic product, with $\bar{\q}_i$ as the quaternionic conjugate of $\q_i$, 
{\it i.e.},
$\tilde{\cal Q}$ is the quaternionic polynomial ring generated by the quaternionic variables $\q_i$.
\et

Similarly, we refer to the ring $\tilde{\cal Q}$ defined in Theorem \ref{thm:2} as the 
{\it quaternionic polynomial ring} generated by quaternionic variables $\q_i$, and call the defining ideal of 
$\tilde{\cal Q}$ the {\it syzygy ideal}, and denote it by 
$\tilde{\cal I}[\q_1, \q_2, \ldots, \q_n]$. 

\bp \label{prop:2}
For all $k>1$, let $\Q=\p_{i_1}\p_{i_2}\cdots \p_{i_k}$. Then
$\p_j[\Q]-[\Q]\p_j$ is in the syzygy ideal $\tilde{\cal I}[\q_1, \q_2, \ldots, \q_n]$.
$[\p_j\Q]-[\Q\p_j]$ is also in the syzygy ideal, and
$A(\Q)$ is in the 3D $\tilde{\mathbb Q}$-vector space
$\tilde{\cal V}^3$.
\ep

\section{Gr\"{o}bner base and normal form: multilinear case for vectors}
\setcounter{equation}{0}

Let $\v_1, \ldots, \v_n$ be vector variables.
For two quaternionic monomials $h_1, h_2$ in the variables, 
$h_1$ is said to be {\it reduced} with respect to $h_2$, if 
$h_1$ is not a {\it multiplier} of $h_2$, or $h_2$ is not a {\it factor} of $h_1$,
{\it i.e.}, there do not exist monomials $l,r$, including rational numbers, such that 
$h_1=lh_2r$. For two quaternionic polynomials $f$ and $g$, $f$ is said to be {\it reduced} with respect to $g$,
if the leading term of $f$ is reduced with respect to that of $g$. The term ``{\it non-reduced}" means the opposite.


Let $\{f_1, f_2, \ldots, f_k\}$ be a set of quaternionic polynomials. A set of quaternionic polynomials
$\{g_1, g_2, \ldots, g_m\}$ is said to be a {\it reduced Gr\"obner base} of the ideal 
${\cal I}:=\langle f_1, f_2, \ldots, f_k \rangle$ generated by the $f_i$ in the quaternionic polynomial ring, if 
(1) $\langle g_1, \ldots, g_m \rangle={\cal I}$, (2) the leading term of any
element in $\cal I$ is a multiplier of the leading term of some $g_i$, (3) the $g_i$ are pairwise reduced
with respect to each other.

The {\it reduction} of a polynomial $f$ with respect to a reduced Gr\"obner base
$g_1, g_2, \ldots, g_m$ is the repetitive procedure of dividing the first non-reduced term $L$ of $f$ by some 
$g_i$ whose leading term is a factor of $L$, then updating $f$ by replacing
$L$ with its remainder, until all terms of $f$ are reduced. 
The result is called the {\it normal form} of $f$ with respect to
the Gr\"obner base. 

By the classical theory of non-commutative Gr\"obner base \cite{mora}, a reduced Gr\"obner base 
always exists for the syzygy ideal of quaternionic polynomials. However, the Gr\"obner base can be computed 
by computers only for small values of $n$ (number of variables) and $d$ (degree of polynomials in the 
Gr\"obner base). Our goal is to find the general form of 
the Gr\"obner base for arbitrary $n$ which can even be a symbol. We first consider the multilinear case. 
 
In a {\it multilinear monomial}, every variable
occurs at most once. A {\it multilinear polynomial} is the sum of multilinear monomials with rational coefficients,
where for every variable, its degree is the same in different monomials. The {\it multilinear product}
is defined for two multilinear polynomials only when they do not have any common variable. 
The {\it multilinear addition} is defined only for multilinear polynomials having the same set of variables.
All multilinear polynomials together with the multilinear product and multilinear addition form an algebraic system,
called the {\it multilinear ring}. All monic multilinear monomials together with the multilinear product form another
algebraic system, called the {\it multilinear monoid}.

An {\it ideal} of a multilinear ring or multilinear monoid is defined just as in the usual tensor-product case.
The concepts of {\it reduced Gr\"obner base}, {\it reduction} and {\it normal form} have similar definitions,
the theory and algorithm of non-commutative Gr\"obner base in \cite{mora} are still valid.

For example, for $m\leq n$ vector variables $\v_{i_1},\ldots, \v_{i_m}$, the {\it multilinear quaternionic 
polynomial ring} 
${\cal Q}^M[\v_{i_1},\ldots, \v_{i_m}]$ is composed of all 
multilinear quaternionic 
polynomials whose variables are among the $\v_{i_k}$, including polynomials of degree 0.
The syzygy ideal ${\cal I}^M[\v_{i_1},$ $\ldots, \v_{i_m}]$ of the ring is generated by 
elements of the form $[\v_{i_a}\v_{i_b}]\v_{i_c}-\v_{i_c}[\v_{i_a}\v_{i_b}]$ 
or $[\v_{i_a}\v_{i_b}\v_{i_c}]\v_{i_d}-\v_{i_d}[\v_{i_a}\v_{i_b}\v_{i_c}]$ for all
$i_a\neq i_b\neq i_c\neq i_d$.
The Gr\"obner base of the syzygy ideal is denoted by ${\cal G}^M[\v_{i_1},\ldots, \v_{i_m}]$.

When the variables $\v_{i_1},\ldots, \v_{i_m}$ are replaced by another sequence $\v_{j_1},\ldots, \v_{j_m}$,
the syzygy ideal of ${\cal Q}^M[\v_{j_1},\ldots, \v_{j_m}]$, and the Gr\"obner base
of the syzygy ideal are obtained from those of ${\cal Q}^M[\v_{i_1},\ldots, \v_{i_m}]$
simply by a replacement of variables. From this aspect, what is important is the number $m$ of
variables, not the specific sequence of variables. We introduce ${\cal Q}^M_m$, ${\cal I}^M_m$,
and ${\cal G}^M_m$ to denote the multilinear quaternionic ring, its syzygy ideal, and the 
Gr\"obner base of the syzygy ideal, in $m$ {\it unspecified vector variables}.

Notice that if ${\cal G}^M_m$ is reduced, then it is a subset of a reduce Gr\"obner base of
${\cal I}^M_n$, as any polynomial of degree $\leq m$ in ${\cal I}^M_n$ is also in ${\cal I}^M_m$.

\bt \label{main:1}
Let 
${\cal I}^M[\v_1, \ldots, \v_n]$ be the syzygy ideal of the multilinear ring 
${\cal Q}^M[\v_1, \ldots, \v_n]$ of 
multilinear quaternionic polynomials in $n$ vector variables $\v_1\prec \v_2\prec \ldots \prec \v_n$. 

(1) [Gr\"obner base] The following are a reduced Gr\"obner base of ${\cal I}^M[\v_1, \ldots, \v_n]$:
for all $1\leq i_1<i_2<\ldots<i_m\leq n$,
\bu
\item[$\rm G3$:] 
${[}\v_{i_3}\v_{i_2}\v_{i_1}]-[\v_{i_1}\v_{i_3}\v_{i_2}]$, and 
${[}\v_{i_3}\v_{i_1}\v_{i_2}]-[\v_{i_2}\v_{i_3}\v_{i_1}]$;

\item[${\rm G}m$:] 
${[}\v_{i_3}\v_{i_2}\V\v_{i_1}]
-{[}\v_{i_2}\V\v_{i_1}\v_{i_3}]$, where
$\V=\v_{i_4}\v_{i_5}\cdots \v_{i_m}$, and $4\leq m\leq n$.
\eu

\vskip -.63cm
\begin{figure}[htbp]
\centering \epsfig{file=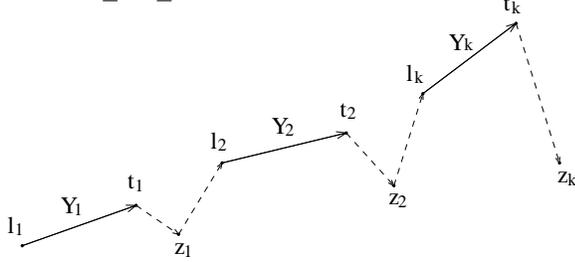, height=3.5cm}
\caption{Double-ascending structure of a multilinear quaternionic monomial in normal form, where only the subscripts 
of vector variables are shown.}
\label{normal}
\end{figure}

(2) [Normal form, see Figure \ref{normal}] In a normal form, every term is
up to coefficient 
of the form $\V_{Y_1}\v_{z_1}\V_{Y_2}\v_{z_2}\cdots \V_{Y_k}\v_{z_k}$ or 
$\V_{Y_1}\v_{z_1}\cdots \V_{Y_{k}}\v_{z_{k}}\V_{Y_{k+1}}$, where \\
(i) $k\geq 0$, \\
(ii) $\v_{z_1}\v_{z_2}\cdots \v_{z_k}$ 
is ascending, \\
(iii) every 
$\V_{Y_i}$ is an ascending monomial of length $\geq 1$, \\
(iv) 
$\V_{Y_1}\V_{Y_2}\cdots \V_{Y_k}$ (or $\V_{Y_1}\V_{Y_2}\cdots \V_{Y_{k+1}}$ if $\V_{Y_{k+1}}$ occurs) is ascending, \\
(v) for every $i\leq k$, if $\v_{t_i}$ is the trailing variable of 
monomial $\V_{Y_i}$, then $\v_{z_i}\prec \v_{t_i}$.
\et

{\it Remark}.\ In a normal form, 
as shown in Figure \ref{normal}, although $\v_{z_i}$ is lower than the trailing variable $\v_{t_i}$ of $\V_{Y_i}$, 
it is allowed that the order of $\v_{z_i}$ be at anywhere from above $\v_{z_{i-1}}$ up to below $\v_{t_i}$, and
both $\v_{z_i}\prec\v_y$ and $\v_{z_i}\succ \v_y$ are possible for any variable $\v_y\neq \v_{t_i}$ in 
$\V_{Y_j}$ where $j\leq i$.

Before proving the theorem, we introduce the following convenient term: let $l_f, l_g$ be the leading terms of
polynomials $f,g$ respectively, if $l_f$ equals $l_g$ up to coefficient, then $f,g$ are said to have
{\it leading like terms}.

\vskip .2cm
\bproof Observe that in the sequence G3, G$m$ for $m\geq 4$, the leading terms are pairwise reduced
with respect to each other. The proof of the theorem is by induction, and is composed of many steps:

{\it Step 1}. G3, G$m$ for $m>3$ generate the ideal ${\cal I}^M_n$. The proof is the following.
By Proposition \ref{prop:1}, G3, G$m$ are in the ideal. It is easy to show that for all $i\neq j\neq k\neq l$,
$[\v_i\v_j\v_k]-[\v_k\v_i\v_j]$ is reduced to zero by G3, and 
$[\v_i\v_j\v_k\v_l]-[\v_l\v_i\v_j\v_k]$ is reduced to G4 up to sign
by G3. So all the generators of ${\cal I}^M_n$ are in the ideal generated by G3, G4.

{\it Step 2}. The theorem obviously holds for $n=3$. For $n=4$, in Step 1 it is shown that for any 
$f\in {\cal I}^M_4$, 
\be
f=\sum_\alpha \lambda_\alpha \v_\alpha \hbox{\rm G3} 
+ \sum_\beta  \lambda_\beta \hbox{\rm G3}\  \v_\beta+\lambda \hbox{\rm G4}, 
\label{form:1}
\ee
where G3 stands for a 
polynomial of type G3, and the $\lambda$'s are rational numbers. Notice that the leading variable of G3 is always
$\v_3$ or $\v_4$. 

Assume that the leading term of $f$ is reduced with respect to G3, G4. Further
assume that before cancelling like terms on the right side of 
(\ref{form:1}), the leading term is $\tilde{L}=\v_{i_1}\v_{i_2}\v_{i_3}\v_{i_4}$ up to coefficient. 
Then $\tilde{L}$ is not reduced with respect to G3, G4, so it must
be cancelled by other leading terms of the addends on the right side of (\ref{form:1}).
The following observations on (\ref{form:1}) are obvious: 

(i) If $\tilde{L}$ comes from $\lambda \rG4$, {\it i.e.}, equals $\v_3\v_2\v_4\v_1$, 
then it is not the leading like term of either $\lambda_\alpha \v_\alpha \rG3$ or $\lambda_\beta \rG3\, \v_\beta$. 
So $\tilde{L}$ cannot come from $\lambda \rG4$.

(ii) If $\tilde{L}$ comes from some $\lambda_i \v_i \rG3$, then it is not the leading like term of any other 
$\lambda_j \v_j \rG3$, so it has to be the leading like term of a unique $\lambda_k \rG3\, \v_k$. Conversely, 
if $\tilde{L}$ comes from some $\lambda_k \rG3\, \v_k$, then it has to be the leading like term of a unique 
$\lambda_i \v_i \rG3$.

Based on the above two observations, we get that in $\tilde{L}=\v_{i_1}\v_{i_2}\v_{i_3}\v_{i_4}$,
$\v_{i_1}\succ \v_{i_2}\succ \v_{i_3}$ and $\v_{i_2}\succ\v_{i_4}$. So either
$\tilde{L}=\v_4\v_3\v_2\v_1$ or $\tilde{L}=\v_4\v_3\v_1\v_2$. 

{\it Step 3}. For $n=4$, we have the following conclusion:
\be
\v_4 \rG3= \sum_{i=1}^3 \lambda_i \v_i \hbox{\rm G3} 
+ \sum_{j=1}^4  \lambda_j \hbox{\rm G3}\  \v_j.
\label{prove:1}
\ee 
To prove this we consider only $h=\v_4(\v_3\v_2\v_1-\v_1\v_2\v_3-\v_1\v_3\v_2+\v_2\v_3\v_1)$. If
(\ref{prove:1}) holds for $\v_4 {\rm G3}=h$, then by interchanging subscripts 1,2, the right side of
(\ref{prove:1}) keeps the same form. 
\[\ba{lcl}
h &\hskip -.45cm=& \hskip -.15cm \phantom{-}
\underline{\v_4\v_3\v_2}\v_1-\underline{\v_4\v_1\v_2}\v_3-
\underline{\v_4\v_1\v_3}\v_2+\underline{\v_4\v_2\v_3}\v_1 \\

&\hskip -.45cm\stackrel{{\rm \scriptstyle G}3\,\v_j}{=} & \hskip -.15cm\phantom{-}
\v_2\v_3\v_4\v_1-\v_2\v_1\v_4\v_3-\v_2\underline{\v_4\v_1\v_3}+\v_1\underline{\v_4\v_2\v_3}\\

&&\hskip -.15cm
-\v_3\underline{\v_4\v_1\v_2}-\v_3\v_1\v_4\v_2+\v_1\underline{\v_4\v_3\v_2}+\v_3\v_2\v_4\v_1\\

&\hskip -.45cm\stackrel{\v_i {\rm \scriptstyle G}3,i<4}{=} & \hskip -.15cm
-\v_2\v_3\v_1\v_4+\v_1\v_3\v_2\v_4+\v_1\v_2\v_3\v_4-\underline{\v_3\v_2\v_1}\v_4
\\

&\hskip -.45cm\stackrel{{\rm \scriptstyle G}3\,\v_4}{=} & \hskip -.15cm
0.
\ea
\]
The underlines give hint on the category of
Gr\"obner base elements to be used in reduction. 

{\it Step 4}. By (\ref{prove:1}), (\ref{form:1}) can be rewritten as
\be
f=\sum_{i<4} \lambda_i \v_i \hbox{\rm G3} 
+ \sum_j  \lambda_j \hbox{\rm G3}\  \v_j+\lambda \hbox{\rm G4}.\label{step4:f}
\ee
By the analysis in Step 2, the leading term on the right side of (\ref{step4:f}) before cancelling like terms
remains the leading term after cancelling like terms, violating the assumption that
the leading term of $f$ is reduced with respect to G3, G4. This proves that G3, G4 form a
Gr\"{o}bner base of ${\cal I}^M_4$.

{\it Step 5}. Let $T=\v_{i_1}\v_{i_2}\v_{i_3}\v_{i_4}$ be a term in the normal form of a degree-4
polynomial. By reduction with G3, if $i_j>i_{j+1}$ then $i_{j+2}>i_j$ is mandatory. So the
sequence of subscripts must be one of $Y_1$, $Y_1z_1$, $Y_1z_1Y_2$, and $Y_1z_1Y_2z_2$, where
$z_i$ is a single subscript, each $Y_i$ is an ascending sequence of subscripts, all the $Y_i$
form a big ascending sequence of subscripts, and $z_i$ is smaller than the trailing
subscript of $Y_i$. If $z_1>z_2$, then it must be that $T=\v_3\v_2\v_4\v_1$, which is not reduced
with respect to G4. So $z_1<z_2$. This proves the theorem for $n=4$.

{\it Step 6}. Assume that the theorem holds for all $n<q$, where $q\geq 5$. When $n=q$, 
let $f\in {\cal I}^M_q$ and assume that its leading term is reduced with respect to 
G3, G$m$ for $m$ up to $q$. By induction hypothesis, the degree of $f$ must be $q$.
We need to prove that $f$ does not exist.

Since the generators of ${\cal I}^M_q$ all have degree $<q$, it must be that
\be
f=\sum_\alpha \lambda_\alpha \v_\alpha f_\alpha
+ \sum_\beta  \lambda_\beta f_\beta \v_\beta, 
\label{form:2}
\ee
where
the $\lambda$'s are coefficients, and $f_\alpha, f_\beta\in {\cal I}^M_{q-1}$ are of degree $q-1$.
We do reduction to $f_\alpha$ with respect to the Gr\"{o}bner base G3, G$m$ for 
$m$ up to $q-1$ in the following special way: to every term $T_i$ of $f_\alpha$, let 
$T_i=\mu_i \V_i\v_{{\alpha}_i}$ where $\V_i$ is a monomial of degree $q-2\geq 3$ and $\mu_i$ is the coefficient,
do reduction to monomial $\V_i$ by Gr\"{o}bner base G3, G$m$ for $m$ up to $q-2$. After the reduction 
to every $\V_i$, 
$f_\alpha=\sum_j  \lambda_j g_j \v_j+\tilde{f}_\alpha$, where $g_j\in {\cal I}^M_{q-2}$, and
$\tilde{f}_\alpha\in {\cal I}^M_{q-1}$. Every term of $\tilde{f}_\alpha$ after removing its trailing variable
is in normal form. 

Now do reduction to $\tilde{f}_\alpha$ by G3, G$m$ for $m$ up to $q-1$. Then
$
\tilde{f}_\alpha=\sum_{m<q}(\sum_{\alpha_m} \lambda_{\alpha_m} a_{\alpha_m} \rG m
+\lambda_{\beta_m} b_{\beta_m} \rG m\ c_{\beta_m}
+\sum_{\gamma_m} \lambda_{\gamma_m} \rG m\ d_{\alpha_m})
$,
where $a,b,c,d$ are monomials of degree $\geq 1$. 
The leading term of $\tilde{f}_\alpha$ must be the leading term of a unique 
$\lambda_{\alpha_m} a_{\alpha_m} \rG m$, because
it can only be non-reduced with respect to such a polynomial. Now substitute the above expressions of 
$\tilde{f}_\alpha$, $f_\alpha$ into (\ref{form:2}), and regroup the terms. We finally get
\be
f=s+t, \hbox{  where  }
s=\sum_{m<q}\sum_{\alpha_m} \lambda_{\alpha_m} \V_{\alpha_m} \rG m, \ \,
t=\sum_\beta  \lambda_\beta h_\beta \v_\beta, 
\label{form:3}
\ee
and where 
$\V_{\alpha_m}$ is a monomial of length $q-m$, 
$h_\beta\in {\cal I}^M_{q-1}$, with the following properties:

(a) If the leading term $T_s$ of an addend
$\lambda_{\alpha_m} \V_{\alpha_m} \rG m$ is of the highest order among all leading terms
in the addends of $s$,
then $T_s$ must be the leading term of $s$, and has no leading like term in any other addend of  
$s$.

(b) Furthermore, when $\V_{\alpha_m}=\v_{i_1}\v_{i_2}\cdots \v_{i_{q-m}}$ and 
the leading term of $\rG m$ is $\v_{j_1}\v_{j_2}\cdots \v_{j_{m}}$, then the degree-$(q-2)$ factor
$\v_{i_2}\cdots \v_{i_{q-m}}\v_{j_1}\v_{j_2}\cdots \v_{j_{m-1}}$ of $T_s$ is in normal form.

{\it Step 7}. In (\ref{form:3}), 
the leading term $T_s$ of $s$ comes from a unique $\lambda_{\alpha_m} \V_{\alpha_m} \rG m$ of $s$. 
We will establish the following conclusion: if $T_s$ has any leading like term in $t$, then 
\be
\V_{\alpha_m} \rG m=\sum_{j\leq q}\sum_{*} \lambda_{\psi_j} \V_{\psi_j} \rG j
+\sum_\phi  \lambda_\phi w_\phi \v_\phi,
\label{establish}
\ee
where the $*$-summation is over all $\V_{\psi_j} \rG j$ whose leading term has lower order than 
that of $\V_{\alpha_m} \rG m$. Once (\ref{establish}) is established, then we can rewrite $f$ as
$\lambda \rG q+\tilde{s}+\tilde{t}$, where 
\be
\tilde{s}=\sum_{j<q}\sum_{\gamma_j} \lambda_{\gamma_j} \V_{\gamma_j} \rG j, \ \
\tilde{t}=\sum_\zeta  \lambda_\zeta \tilde{h}_\zeta \v_\zeta, 
\label{form:4}
\ee
such that either $\tilde{s}=0$, or the leading term of $\tilde{s}$ does not have 
leading like term in $\tilde{t}$.

The two observations in Step 2 are also valid in the current situation. The leading term of 
$\tilde{s}$ can only be cancelled by a term of $\tilde{t}$, and conversely, the
leading term of $\tilde{t}$ can only be cancelled by a term of $\tilde{s}$.
Furthermore, the leading terms of $\tilde{s}$ and $\tilde{t}$ are both non-reduced
with respect to G3, G$j$ for $j$ up to $q-1$. Since  
the leading term of $f$ is reduced with respect to G3, G$j$ for $j$ up to $q-1$, it must have higher order
than any term of $\tilde{s}$ and $\tilde{t}$; then it must be the leading term of G$q$, 
violating the assumption that it is reduced with respect to G$q$. This proves that
G3, G$j$ for $j$ up to $q$ is a Gr\"{o}bner base once (\ref{establish}) is true. 

{\it Step 8}. Once it is proved that G3, G$j$ for $j$ up to $q$ is a Gr\"{o}bner base, then for
a polynomial $f$ of degree $q$, by reduction with respect to 
G3, G$j$ for $j$ up to $q-1$, we get that for any term $T=\lambda \v_{l}\V_{q-2}\v_t$ in the result of reduction,
where $\V_{q-2}$ is a monomial of degree $q-2$, the subscripts of
$\V_{q-2}\v_t$ and $\v_{l}\V_{q-2}$ must be both in one of the following normal forms: 
$Y_1$, $Y_1z_1$, $Y_1z_1Y_2$, $Y_1z_1Y_2z_2$, $\ldots$. If the subscripts of
$T$ are not in one of the normal forms, it is easy to show that (1) the subscripts of
$\V_{q-2}\v_t$ must be of the form $Y_1z_1$, and $z_1=t$; (2) if $Y_1=h_1h_2Y$, where $Y$ is of length $q-4\geq 1$, then 
$h_2>l>h_1$; (3) $h_1>t$. The only possibility is $\v_{l}\V_{q-2}\v_t=\v_3\v_2\v_4\v_5\cdots \v_{q}\v_1$, 
the leading term of G$q$.

Consider polynomial G$q$. By reduction with respect to G3, G$j$ for $j$ up to $q-1$, we get G$q
=\v_3\v_2\v_4\v_5\cdots \v_{q}\v_1-\v_3\v_1\v_4\v_5\cdots \v_{q}\v_2+g$, where each term of $g$ has leading variable
$\v_1$ or $\v_2$. The second term $\v_3\v_1\v_4\v_5\cdots \v_{q}\v_2$ and the terms of $g$ are each in a normal form
described by the current theorem. 
Hence, after reduction with respect to the G$q$ in the above form, each term of $f$ becomes normal.

{\it Step 9}. We start to prove (\ref{establish}), which is the main part of the whole proof.
Let $\V_{\alpha_m}=\v_w\V_{q-m-1}$, and let
the leading term of $\rG m$ be $\V^G_{m-1}\v_z$, then up to scale $T_s=\v_w\V_{q-m-1}\V^G_{m-1}\v_z$ is 
the leading term of $s$, and monomial $\V_{q-m-1}\V^G_{m-1}$ is in normal form, {\it i.e.}, its
subscripts are in one of the forms $Y_1z_1\ldots Y_kz_k$ and $Y_1z_1\ldots Y_kz_kY_{k+1}$. 
Let $l_i, t_i$ be the leading subscript and trailing subscript of $Y_i$ respectively. Then
$l_i\leq t_i$. Let the middle part of $Y_i$ after removing $l_i, t_i$ be $C_i$.

The tail part of the subscripts $Y_1z_1\ldots Y_kz_kz$ or 
$Y_1z_1\ldots $ $Y_kz_kY_{k+1}z$ must be either the subscripts $y\tilde{z}z$ from the leading term of
G3 for $m=3$, where $y>z$ and $y>\tilde{z}$,
or the subscripts $y\tilde{z}Yz$ from the leading term of G$j$ for $m=j>3$, where 
$z\tilde{z}yY$ is ascending.
So subscripts $Y_1z_1\ldots Y_{k-1}z_{k-1}$ $(l_kC_kt_k)z_kz$ where $t_k>z$, belong to the leading term
of $\V_{Y_1}\v_{z_1}\ldots \V_{Y_{k-1}}\v_{z_{k-1}}\v_{l_k}\V_{C_k}\rG 3$, while subscripts
$Y_1z_1\ldots$ $ Y_{k-1}z_{k-1}(l_kC_kt_k)z_kY_{k+1}z$ where $z_k>z$, belong to the leading term
of $\V_{Y_1}\v_{z_1}$ $\cdots \V_{Y_{k-1}}\v_{z_{k-1}}\v_{l_k}\V_{C_k}\rG j$, where $j$ is the length of $Y_{k+1}$
plus 3.
We get
\be
\V_{\alpha_m}=\v_w\V_{Y_1}\v_{z_1}\cdots \V_{Y_{k-1}}\v_{z_{k-1}}\v_{l_k}\V_{C_k}.
\ee 

If $\V_{\alpha_m}$ is not in normal form,
then by reduction with G3, G$j$ for $j$ up to $q-m$,
the order of $\V_{\alpha_m}$ is decreased, and (\ref{establish}) is trivially
satisfied. Below we assume that $\V_{\alpha_m}$ is normal. 

Since $T_s$ has leading like term in $t$,
$\v_w\V_{\alpha_m}\v_{t_k}\v_{z_k}$ for $t_k>z$, and
$\v_w\V_{\alpha_m}\v_{t_k}\v_{z_k}\V_{Y_{k+1}}$ for $z_k>z$, must be non-reduced with respect to 
G3, G$j$ for $j$ up to $q-1$. It is easy to deduce $k=1$, and when any of the two monomials
is non-reduced with respect to G3, the length of $Y_1$ is less than 3. So up to scale,
\be
T_s=\v_w\V_{Y_1}\v_{z_1}\v_z\ \hbox{   or   }\ \v_w\V_{Y_1}\v_{z_1}\V_{Y_2}\v_z.
\ee
Only the head part and tail part
of $T_s$ are non-reduced with respect to some G$i$, G$j$ respectively, where $i,j$ are between 3 and $q-1$. 
There are the following cases:

(I) $i=j=3$. The subscripts of $T_s$ are either $wt_1z_1z$ or 
$wl_1t_1z_1z$, where $w>t_1$ and $t_1>z$. The first case is considered in Step 3 and is conclusion (\ref{prove:1}).
The second case includes six different subscripts of $T_s$: 51423, 51432, 52413, 52431, 53412, 53421. 

(II) $i=3$, $j>3$. The subscripts of $T_s$ are either $wt_1z_1Y_2z$ or $wl_1t_1z_1Y_2z$,
where $w>t_1$ and $z_1>z$. 

(II.1) In the former case, the subscripts of $T_s$ are 
$w32Y1$, where $w>3$, and $Y>3$ is an 
ascending sequence of length $q-4>0$.

(II.2)
In the latter case, $t_1=4$, and $t_1,z_1,z$ is a permutation of 1,2,3. So
$T_s$ has the following possible subscripts: $w143Y2$,
$w243Y1$, and $w342Y1$, where $w>4$, and $Y>4$ is an ascending sequence of length $q-5>0$.

(III) $i>3$, $j=3$. The subscripts of $T_s$ are $wl_1C_1t_1z_1z$, where 
$C_1$ is an ascending sequence of length $q-5>0$, 
$w>l_1>z_1$ but $w<C_1$, and $t_1>z$. 
Obviously $t_1=q$, but $z$ can take any value from 1 to $q-1$. So
$T_s$ has the following possible subscripts: $4356\ldots q21$, $4356\ldots q12$, $4256\ldots q13$, 
$3256\ldots q14$, and $3245\ldots \check{k}\ldots q1k$ for $5\leq k\leq q-1$. 
There are the following subcases, considering both the symmetries of the subscripts and the difference
of the categories of Gr\"obner base elements to be used in reduction:

(III.1) The subscripts of $T_s$ are 
$43521$, $43512$, or $42513$.

(III.2) The subscripts of $T_s$ are 
$43Y(q-1)q21$, $43Y(q-1)q12$, or $42Y(q-1)q13$, where $Y$ is an ascending sequence of length $q-6\geq 0$, 
and $4<Y<q-1$.

(III.3)
The subscripts of $T_s$ are $32514$.

(III.4)
The subscripts of $T_s$ are 
$32Yq1(q-1)$, where $Y$ is an ascending sequence of length $q-5\geq 0$, 
and $3<Y<q-1$. 

(III.5)
The subscripts of $T_s$ are 
$32Y(q-1)q14$, where $Y$ is an ascending sequence of length $q-6\geq 0$, and $4<Y<q-1$.

(III.6)
The subscripts of $T_s$ are 
$32YZ(q-1)q1k$, where $4<k<q-1$, $Y$ is an ascending sequence of length $u>0$, 
$Z$ is an ascending sequence of length
$q-u-6\geq 0$, and $3YkZ(q-1)$ is an ascending sequence. 

(IV) $i>3$, $j>3$. The subscripts of $T_s$ are $wl_1C_1t_1z_1Y_2z$, where 
$C_1, Y_2$ are two ascending sequences of length $>0$, 
$w>l_1>z_1>z$ but $w<C_1$. 
So $zz_1l_1w=1234$, and 
the subscripts of $T_s$ are $43Yk2Z1$, where $Y$ is an ascending sequence of length $u\geq 0$, 
$Z$ is an ascending sequence of length $q-u-5>0$, and $4YkZ$ is an ascending sequence.

There is no space to present the analyses of all the cases and subcases. We select to present the details
for the following representative cases: (I), (II.1), (III.2) to (III.4). 


{\it Step 10}. The polynomials on the left side of the following equalities
are in ${\cal I}^M_{j+k}$: for monomials $\V_j, \V_k$ of degree $j,k$ respectively,
\be\left\{\ba{rll}
(\V_k+(-1)^k\V_k^\dagger)\V_j
-\V_j(\V_k+(-1)^k\V_k^\dagger) &=& 0,\\
\V_j\V_k-\V_k\V_j
-(-1)^{j+k}(\V_j^\dagger\V_k^\dagger-\V_k^\dagger\V_j^\dagger) &=& 0.
\ea\right.
\label{reduction}
\ee
As long as $j+k<q$, and by applying the above equalities to the leading term of a polynomial $f$,
the order of $f$ is decreased, this transformation of $f$ can always be replaced by reduction with 
respect to G3, G$i$ for $i$ up to $j+k$. 

Let both $\v_{i_1},\v_{i_2}\prec \v_{i_3}$ and $\v_{i_1},\v_{i_2}\prec \V_A$, where $A$ is a sequence of 
subscripts of length $a<q-3$. 
The following order-reducing transformation is also induced by G3, G$i$ for $i$ up to $a+3$:
\be\ba{lll}
\v_{i_3}\V_A\v_{i_1}\v_{i_2} &=& \v_{i_2}(\V_A\v_{i_1}\v_{i_3}+(-1)^a \v_{i_3}\v_{i_1}\V_A^\dagger)\\
&& \hfill -(-1)^a \v_{i_1}\V_A^\dagger\v_{i_3}\v_{i_2}.
\ea
\label{reduction:2}
\ee
The proof of (\ref{reduction:2}) is simple:
\[\ba{ll}
&\hskip -.25cm  \v_{i_3}\V_A\v_{i_1}\v_{i_2} \\

=&\hskip -.25cm \v_{i_3}(\V_A\v_{i_1}\v_{i_2}+(-1)^a \v_{i_2}\v_{i_1}\V_A^\dagger)
-(-1)^a\underline{\v_{i_3}\v_{i_2}\v_{i_1}}\V_A^\dagger
\\

=&\hskip -.25cm  (\v_{i_2}\V_A\v_{i_1}+(-1)^a \v_{i_1}\V_A^\dagger\v_{i_2})\v_{i_3}
+(-1)^a\v_{i_2}\v_{i_3}\v_{i_1}\V_A^\dagger\\

&\hfill-(-1)^a\v_{i_1}\V_A^\dagger(\v_{i_3}\v_{i_2}+\v_{i_2}\v_{i_3})
\\

=&\hskip -.25cm  \v_{i_2}(\V_A\v_{i_1}\v_{i_3}
+(-1)^a\v_{i_3}\v_{i_1}\V_A^\dagger)
-(-1)^a\v_{i_1}\V_A^\dagger\v_{i_3}\v_{i_2}.
\ea
\]
In later steps, (\ref{reduction}) and (\ref{reduction:2}) are used directly in reduction.

{\it Step 11}. Consider case (I) of Step 9.
Let
$h=\v_5\v_1(\v_4\v_2\v_3$ $+\v_2\v_4\v_3-\v_3\v_4\v_2-\v_3\v_2\v_4)$. Then
\[\hskip -.2cm
\ba{lcl}
h &\hskip -.45cm=& \hskip -.45cm 
\underline{\v_5\v_1\v_4}\v_2\v_3+\underline{\v_5\v_1\v_2}\v_4\v_3-
\underline{\v_5\v_1\v_3}(\v_2\v_4+\v_4\v_2) \\

&\hskip -.45cm \stackrel{{\cal \scriptstyle I}_4^M\v_i,i<5}{=} & \hskip -.25cm
-\underline{\v_1\v_5(\v_4\v_2+\v_2\v_4)\v_3}
+\v_4(\v_5\v_1+\v_1\v_5)\v_2\v_3
\\

&& \hskip -.25cm
+\v_2(\v_5\v_1+\v_1\v_5)\v_4\v_3
+\underline{\v_1\v_5\v_3(\v_2\v_4+\v_4\v_2)}
\\

&& \hskip -.25cm
-\v_3(\v_5\v_1+\v_1\v_5)(\v_2\v_4+\v_4\v_2)
\\

&\hskip -.45cm \stackrel{\v_j{\cal\scriptstyle I}_4^M,j<5}{=} & \hskip -.25cm
(\v_4\v_2\v_3+\v_2\v_4\v_3
-\v_3\v_2\v_4-\v_3\v_4\v_2)\\

&&
\hfill
(\v_5\v_1+\v_1\v_5)
\\

&\hskip -.45cm\stackrel{{\cal \scriptstyle I}_4^M\v_k}{=} & 
\hskip -.25cm
0.
\ea
\]
When the subscripts 1,2,3 undergo a permutation, the polynomial equalities used in the reduction
remain in the same category:
${\cal I}_4^M\v_i$ for $i<5$; $\v_j{\cal I}_4^M$ for $j<4$; ${\cal I}_4^M\v_k$. So the above
reduction is valid for all the six subcases of Case (I).

{\it Step 12}. Case (II.1) of Step 9. 
Let 
$h=\v_w(\v_3\v_2
\V_Y\v_1+(-1)^q\v_3\v_1\V_Y^\dagger\v_2-(-1)^q\v_1\V_Y^\dagger\v_2\v_3
-\v_2\V_Y\v_1\v_3)$, where $w>3$, and $Y>3$ is an 
ascending sequence of length $q-4>0$. Then
\[\ba{lcl}
h &\hskip -.25cm \stackrel{{\cal\scriptstyle I}_{q-1}^M\v_i, i<3}{=} & \hskip -.25cm
-\v_3\v_w(\v_2\V_Y\v_1+(-1)^q \v_1\V_Y^\dagger\v_2)\\

&& \hskip -.25cm+\v_2(\v_w\v_3+\v_3\v_w)\V_Y\v_1
\\

&& \hskip -.25cm
+(-1)^q \v_1(\v_w\v_3+\v_3\v_w)\V_Y^\dagger\v_2\\

&& \hskip -.25cm
-\v_w(\v_2\V_Y\v_1+(-1)^q \v_1\V_Y^\dagger\v_2)\v_3
\\

&\hskip -.25cm \stackrel{\v_j{\cal\scriptstyle I}_{q-1}^M, j<4}{=} & \hskip -.25cm
-\v_3(\v_2\V_Y\v_1+(-1)^q \v_1\V_Y^\dagger\v_2)\v_w\\

&& \hskip -.25cm
+(\v_2\V_Y\v_1+(-1)^q \v_1\V_Y^\dagger\v_2)(\v_w\v_3+\v_3\v_w)
\\

&& \hskip -.25cm
-\v_w(\v_2\V_Y\v_1+(-1)^q \v_1\V_Y^\dagger\v_2)\v_3
\\

&\hskip -.25cm \stackrel{{\cal \scriptstyle I}_{q-1}^M\v_j}{=} & \hskip -.25cm
\hskip -.25cm
0.
\ea
\]


\def\localgap{\hskip -1cm}
\def\equalgap{\hskip -.65cm}

{\it Step 13}. Case (III.2) of Step 9. 
Let 
$h=\v_4\v_3\V_Y\v_{q-1}(\v_q\v_2$ $\v_1-\v_1\v_2\v_q
-\v_1\v_q\v_2+\v_2\v_q\v_1)$, 
where $Y$ is an ascending sequence of length $q-6\geq 0$, 
and $4<Y<q-1$. Then
\[\hskip -.14cm
\ba{lcl}
h &\equalgap = & \localgap
\underline{\v_4\v_3\V_Y\v_{q-1}\v_q\v_2}\v_1-\v_4\v_3\V_Y\v_{q-1}\underline{\v_1\v_2}\v_q\\

&& \localgap
-\v_4\v_3\V_Y\v_{q-1}\v_1\v_q\v_2+\v_4\v_3\V_Y\underline{\v_{q-1}\v_2\v_q\v_1}\\


&\equalgap \stackrel{\stackrel{\V_{43Y}{\scriptstyle \rm G}4,}{{\cal\scriptstyle I}_{q-1}^M\v_i}}{=} & \localgap
\hskip .1cm
(\v_3\V_Y\v_{q-1}\v_q\v_2+(-1)^q\v_2\v_q\v_{q-1}\V_Y^\dagger\v_3)\v_4\v_{1}\\

&& \localgap
-(-1)^q \v_4\v_2\underline{\v_q\v_{q-1}\V_Y^\dagger\v_3\v_1}\\

&& \localgap
-\v_2(\v_3\V_Y\v_{q-1}\v_1\v_4+(-1)^q\v_4\v_1\v_{q-1}\V_Y^\dagger\v_3)\v_{q}\\

&& \localgap
+(-1)^q \v_1\v_{q-1}\V_Y^\dagger\v_3\v_4\v_2\v_{q}\\

&& \localgap
+\v_4\v_3\V_Y\underline{(\v_2\v_q\v_1-\v_1\v_q\v_2)}\v_{q-1}\\

&\equalgap \stackrel{\stackrel{\V_{42}{\cal\scriptstyle I}_{q-2}^M,{\cal\scriptstyle I}_{q-1}^M\v_i,}
{\v_2{\cal\scriptstyle I}_{q-1}^M}}{=} & \localgap
\hskip .2cm
(\v_3\V_Y\v_{q-1}\v_q\v_2+(-1)^q\v_2\v_q\v_{q-1}\V_Y^\dagger\v_3)\v_4\v_{1}\\

&& \localgap
-\underline{\v_4\v_2(\v_3}\V_Y\v_{q-1}\v_q\v_1+(-1)^q\v_1\v_q\v_{q-1}\V_Y^\dagger\v_3)\\

&& \localgap
+\v_4\v_2\v_1\v_3\V_Y\v_{q-1}\v_q\\

&&\localgap
-\v_2\v_q(\v_1\v_4\v_3\V_Y\v_{q-1}+(-1)^q\v_{q-1}\V_Y^\dagger\v_3\v_4\v_1)\\

&& \localgap
+(-1)^q \v_1\v_{q-1}\V_Y^\dagger\v_3\v_4\v_2\v_{q}\\

&&\localgap
+(\v_2\v_q\v_1-\v_1\v_q\v_2)\v_4\v_3\V_Y\v_{q-1}\\

&\equalgap \stackrel{{\cal\scriptstyle I}_{q-1}^M\v_i,
\v_1{\cal\scriptstyle I}_{q-1}^M}{=} & \localgap
\Bigstrut
-\v_1(\V_Y\v_{q-1}\v_q\v_2\v_4\v_3-(-1)^q\v_3\v_4\v_2\v_q\v_{q-1}\V_Y^\dagger)\\

&& \localgap
+\v_2\v_4\v_3\V_Y\v_{q-1}\v_q\v_1-\v_3\V_Y\v_{q-1}\v_q\v_4\underline{\v_2\v_1}\\

&& \localgap
+\v_4\v_2\v_1(\v_3\V_Y\v_{q-1}\v_q-(-1)^q\v_q\v_{q-1}\V_Y^\dagger\v_3)\\

&\equalgap \stackrel{\v_3{\cal\scriptstyle I}_{q-1}^M,
\v_1{\cal\scriptstyle I}_{q-1}^M}{=} & \localgap 
\Bigstrut
-\v_1(\V_Y\v_{q-1}\v_q\v_2\v_4\v_3-(-1)^q\v_3\v_4\v_2\v_q\v_{q-1}\V_Y^\dagger)\\

&& \localgap
-\v_1(\V_Y\v_{q-1}\v_q\v_4\v_2\v_3-(-1)^q\v_3\v_2\v_4\v_q\v_{q-1}\V_Y^\dagger)\\

&& \localgap
+\v_2\v_4(\v_3\V_Y\v_{q-1}\v_q-(-1)^q\v_q\v_{q-1}\V_Y^\dagger\v_3)\v_1
\\

&& \localgap
+\underline{\v_4\v_2\v_1}(\v_3\V_Y\v_{q-1}\v_q-(-1)^q\v_q\v_{q-1}\V_Y^\dagger\v_3)
\\

&\equalgap \stackrel{{\cal\scriptstyle I}_{q-1}^M\v_i,
\v_1{\cal\scriptstyle I}_{q-1}^M}{=} & \localgap
\Bigstrut \hskip .6cm
0.
\ea
\]


In the reduction rules, $\V_{43Y}$G4 stands for $\v_4\v_3\V_Y$G4, and so for the rest.
The reductions for the other two subcases of Case (III.2) are similar.

{\it Step 14}. Case (III.3) of Step 9. 
Let $h=\v_3\v_2(\v_5\v_1\v_4-\v_4\v_1\v_5-\v_4\v_5\v_1+\v_1\v_5\v_4)$. Then 
\[\ba{lcl}
h &
\hskip -.55cm \stackrel{{\cal \scriptstyle I}_4^M\v_i,{\scriptstyle \rm G}5}{=} & \hskip -.45cm
\v_3\v_1\underline{\v_5(\v_2\v_4+\v_4\v_2)}+(\v_2\v_5\v_1-\v_1\v_5\v_2)\v_3\v_4\\

&& \hskip -.25cm
-\v_3\v_1\v_4\v_2\v_5-(\v_2\v_4\v_1-\v_1\v_4\v_2)\v_3\v_5
\\

&& \hskip -.25cm
-\v_2\v_4\v_5\v_1\v_3-\v_1\v_5\v_4\v_2\v_3+\underline{\v_3\v_2\v_1}\v_5\v_4
\\

&\hskip -.55cm 
\stackrel{\V_{31}{\scriptstyle \rm G}3,{\cal \scriptstyle I}_4^M\v_4}{=} & \hskip -.25cm
\v_3\v_1\v_2\v_4\v_5-\v_2\v_3\v_1\v_5\v_4\\

&& \hskip -.25cm
-\v_1\v_5\v_2\v_3\v_4+\v_1\v_4\v_2\v_3\v_5-\v_1\v_5\v_4\v_2\v_3
\\

&& \hskip -.25cm
+\v_2\v_5\v_1\v_3\v_4-\v_2\v_4\v_1\v_3\v_5-\v_2\v_4\v_5\v_1\v_3
\\

&& \hskip -.25cm
+\v_1(\v_2\v_3+\v_3\v_2)\v_5\v_4
\\

&\hskip -.55cm \stackrel{{\cal\scriptstyle I}_4^M\v_5,\v_1{\cal\scriptstyle I}_4^M}{=} & \hskip -.25cm
\v_2(\v_1\v_3+\v_3\v_1)\v_4\v_5-\v_2\v_4\v_1\v_3\v_5
\\

&& \hskip -.25cm
+\v_2(\v_5\v_1\v_3-\v_3\v_1\v_5)\v_4-\v_2\v_4\v_5\v_1\v_3
\\

&\hskip -.55cm\stackrel{\v_2{\cal \scriptstyle I}_4^M}{=} & 
\hskip -.25cm
0.
\ea
\]

{\it Step 15}. Case (III.4) of Step 9. 
Let 
$h=\v_3\v_2\V_Y(\v_q\v_1$
$\v_{q-1}-\v_{q-1}\V_1\v_q
-\v_{q-1}\v_q\v_1+\v_1\v_q\v_{q-1})$, where 
$Y$ is an ascending sequence of length $q-5\geq 0$, and $3<Y<q-1$. Then
\[\hskip -.14cm
\ba{lcl}
h &\hskip -.25cm \stackrel{{\scriptstyle \rm G}q,{\cal\scriptstyle I}_{q-1}^M\v_i}{=} & \hskip -.25cm
(\v_2\V_Y\v_q\v_1+(-1)^q\v_1\v_q\V_Y^\dagger\v_2)\v_3\v_{q-1}\\

&& \hskip -.25cm
-(-1)^q \v_3\v_1\underline{\v_q\V_Y^\dagger\v_2}\v_{q-1}\\

&& \hskip -.25cm
-(\v_2\V_Y\v_{q-1}\v_1+(-1)^q\v_1\v_{q-1}\V_Y^\dagger\v_2)\v_3\v_{q}\\

&& \hskip -.25cm
+(-1)^q \v_3\v_1\underline{\v_{q-1}\V_Y^\dagger\v_2}\v_{q}\\

&& \hskip -.25cm
-(\v_2\V_Y\v_{q-1}\v_q\v_1-(-1)^q\v_1\v_q\v_{q-1}\V_Y^\dagger\v_2)\v_3\\

&& \hskip -.25cm
-(-1)^q \v_3\v_1\underline{\v_q\v_{q-1}\V_Y^\dagger\v_2}\\

&& \hskip -.25cm
+(\v_2\V_Y\v_1-(-1)^q\v_1\V_Y^\dagger\v_2)\v_3\v_{q}\v_{q-1}\\

&& \hskip -.25cm
+(-1)^q \v_3\v_1\underline{\V_Y^\dagger\v_2}\v_{q}\v_{q-1}\\

&\hskip -.25cm \stackrel{\stackrel{\V_{31}{\cal\scriptstyle I}_{q-2}^M;}
{\v_i{\cal\scriptstyle I}_{q-1}^M,i<3}}{=} & \hskip -.25cm
\v_3\v_1\V_Y(\v_{q}\v_2\v_{q-1}-\v_{q-1}\v_2\v_{q}-\v_{q}\v_{q-1}\v_2\\

&& \hskip -.25cm
+\v_2\v_{q}\v_{q-1})+\underline{\v_3\v_1\v_2}\V_Y(\v_{q-1}\v_{q}-\v_{q}\v_{q-1})\\

&& \hskip -.25cm
-\v_2(\v_1\v_3+\v_3\v_1)\V_Y(\v_{q-1}\v_{q}-\v_{q}\v_{q-1})\\

&& \hskip -.25cm
+\v_1\v_3\v_2\V_Y(\v_{q-1}\v_{q}-\v_{q}\v_{q-1})\\

&\hskip -.25cm \stackrel{\stackrel{\V_{31}{\cal\scriptstyle I}_{q-2}^M,}
{{\cal\scriptstyle I}_{q-1}^M\v_i}}{=} & \hskip -.25cm
0.
\ea
\]




\section{Gr\"{o}bner base and normal form for vector variables}
\setcounter{equation}{0}

\bt \label{main:2}
Let ${\cal I}[\v_1, \ldots, \v_n]$ be the syzygy ideal of the quaternionic polynomial ring ${\cal Q}[\v_1, \ldots, \v_n]$
in $n$ vector variables $\v_1\prec \v_2\prec \ldots \prec \v_n$. 

(1) [Gr\"obner base] The following are a reduced Gr\"obner base of ${\cal I}[\v_1, \ldots, \v_n]$:
for all $1\leq i_1<i_2<i_3\leq i_4\leq \ldots \leq i_{m-1}<i_m\leq n$, 
\bu
\item[$\rm VG3$:] 
\[\ba{cc}
{[}\v_{i_3}\v_{i_2}\v_{i_1}]-[\v_{i_1}\v_{i_3}\v_{i_2}], &
{[}\v_{i_3}\v_{i_1}\v_{i_2}]-[\v_{i_2}\v_{i_3}\v_{i_1}],\\
\v_{i_2}\v_{i_2}\v_{i_1}-\v_{i_1}\v_{i_2}\v_{i_2}, &
\v_{i_2}\v_{i_1}\v_{i_1}-\v_{i_1}\v_{i_1}\v_{i_2}; \\
\ea
\]

\item[${\rm VG}m$:]
${[}\v_{i_3}\v_{i_2}\V\v_{i_1}]
-{[}\v_{i_2}\V\v_{i_1}\v_{i_3}]$, where
$\V=\v_{i_4}\v_{i_5}\cdots \v_{i_m}$, and $m\geq 4$.
\eu

\begin{figure}[htbp]
\hskip -.34cm
\centering \epsfig{file=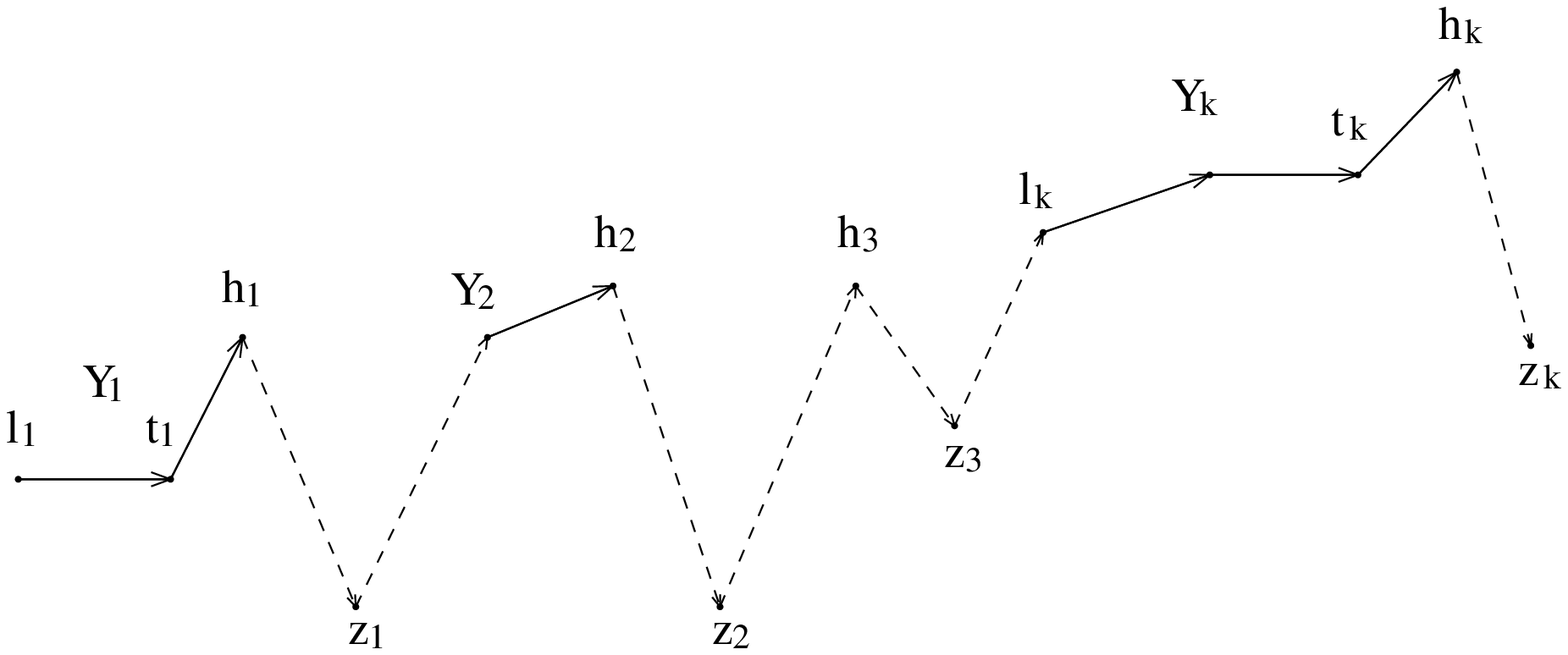, height=3.65cm}
\caption{Double-nondescending structure of a quaternionic monomial in normal form, where only the subscripts 
of vector variables are shown. Each $h_i$ is a peak, while each $z_i$ is a bottom.}
\label{normal:2}
\end{figure}

(2) [Normal form, see Figure \ref{normal:2}] In a normal form, every term is up to coefficient
of the form $\V_{Y_1}\v_{h_1}\v_{z_1}\V_{Y_2}\v_{h_2}\v_{z_2}$ $\cdots \V_{Y_k}\v_{h_k}\v_{z_k}$ or 
$\V_{Y_1}\v_{h_1}\v_{z_1}\cdots \V_{Y_k}\v_{h_k}\v_{z_k}\V_{Y_{k+1}}$,
where \\
(i) $k\geq 0$, \\
(ii) $\v_{z_1}\v_{z_2}\cdots \v_{z_k}$ is non-descending, \\
(iii) $\v_{h_1}\v_{h_2}\cdots \v_{h_k}$ 
is non-descending, \\
(iv) every 
$\V_{Y_i}$ is a non-descending monomial of length $\geq 0$, \\
(v) 
$\V_{Y_1}\v_{h_1}\V_{Y_2}\v_{h_2}\cdots \V_{Y_k}\v_{h_k}$ 
(or $\V_{Y_1}\v_{h_1}\cdots \V_{Y_k}\v_{h_k}\V_{Y_{k+1}}$
if $\V_{Y_{k+1}}$ occurs) is non-descending, 
\\
(vi) for every $i\leq k$, $\v_{h_i}\succ \v_{z_i}$,
\\
(vii) for every $i\leq k$, 
if the length of $\V_{Y_i}$ is nonzero, let
$\v_{t_i}$ be the trailing variable of $\V_{Y_i}$, then
$\v_{t_i}\prec \v_{h_i}$.
\et

\noindent{\bfseries\itshape Sketch of the proof.}\hskip1pc\ignorespaces
Since the generators V2, V3, V4 of ${\cal I}[\v_1, \ldots, \v_n]$ are homogeneous polynomials, 
if $f\in {\cal I}[\v_1, \ldots$, $\v_n]$, then
by writing $f=\sum_\alpha g_\alpha$, where each $g_\alpha$ is a homogeneous polynomial, and different
$g$'s have different multisets of variables, we get that each $g_\alpha\in {\cal I}[\v_1, \ldots, \v_n]$. So we simply 
consider homogeneous polynomials of degree $\leq d$ in $r$ vector variables, where $d\geq r$. Let the multiset
of the variables in such polynomials be ${\cal M}$, whose number of elements is $d$. 

The elements of ${\cal M}$ are specifications of $d$ different vector variables $\cal D$. The definitions of 
multilinear product and multilinear addition of the multilinear ring for $\cal D$ are valid for $\cal M$ after minor revisions.
Let ${\cal I}^M[{\cal M}]$ be the syzygy ideal of the multilinear ring ${\cal Q}^M[{\cal M}]$ of 
multilinear polynomials whose variables by counting multiplicity are in multiset $\cal M$. 
We only need to prove the theorem for ${\cal I}^M[{\cal M}]$. 

Most steps of the proof of Theorem \ref{main:1} are still valid, in particular, (\ref{establish}) is valid, whose proof
follows a similar classification as in Step 9. Although there are more subtle cases to
consider, the reductions for the new situations are much the same as for the cases of Step 9.

\vskip .2cm
This paper is supported partially by NSFC 10871195, \\
60821002/F02, and NCMIS of CAS.

\end{document}